# Extensions to the Guaranteed Service Model for Industrial Applications of Multi-Echelon Inventory Optimization

Victoria G. Achkar[a,b], Braulio B. Brunaud[c], Héctor D. Pérez[d], Rami Musa[c], Carlos A. Méndez[a,b], Ignacio E. Grossmann[d]

[a] *Universidad Nacional del Litoral, Argentina*
[b] *INTEC (UNL – CONICET), Argentina*
[c] *Johnson & Johnson, USA*
[d] *Carnegie Mellon University, USA*

**Abstract**. Multi-echelon inventory optimization (MEIO) plays a key role in a supply chain seeking to achieve specified customer service levels with a minimum capital in inventory. In this work, we propose a generalized MEIO model based on the Guaranteed Service approach to allocate safety stock levels across the network at the lowest holding cost. This model integrates several existing and some novel features that are usually present in pharmaceutical multi-echelon supply chains into a single model: review periods, manufacturing facilities, hybrid nodes (nodes with both internal and external demand), minimum order quantities (MOQ), and different service level performance indicators (fill rate and cycle service levels). We include a polynomial regression to approximate fill rates as a possible target measure to set safety stocks. To improve efficiency, we propose a nonlinear programming model to support decision making, which can be reformulated as a Quadratically Constrained Program (QCP), which leads to order of magnitude reductions in computational time. The performance of the model is evaluated by solving illustrative and real-world cases, and is validated with simulation.

1. INTRODUCTION

On-time fulfilment of customer demand is critical in today's customer-centric supply chains. Achieving this goal depends in great part on the inventory levels and policies that are set along a supply chain. However, efficient inventory control is particularly challenging when customer demand is uncertain and retailers may not know the exact size of an order in advance. Other sources of uncertainty may increase the problem complexity, such as lead time variability. Moreover, the decision at one stage impacts inventory decisions at other stages. The intent of safety stock allocation is to determine an overall strategy for deploying inventory levels across the supply chain in order to buffer it against sources of uncertainty (Graves & Willems, 2003). To overcome these challenges, having a safety stock serves to mitigate the risk of stock-outs in the system. The purpose is to allocate safety stocks to meet customer service levels, while minimizing the total capital tied up in inventory throughout the supply chain, in contrast to single-echelon inventory optimization (SEIO), which seeks to independently minimize cost at each echelon. MEIO has enabled companies to reduce their inventories up to 30% and improve item availability up to 5% by supporting supply chain segmentation, and providing a better balance between lead time, inventory and service under uncertainty (Gartner, 2016). From an optimization perspective, making decisions about inventory in multi-echelon systems is a challenging task because the objective functions usually involve non-linearities, and decision variables affect more than one echelon.

MEIO approaches have been studied in the literature for allocating safety stock in supply chains. De Kok et al. (2018) present a general typology and review stochastic MEIO models in which they classify the extensive research on multi-echelon inventory management under model assumptions, research goals, and different applied methodologies. They state that multi-echelon inventory systems are still a very active area of research because of their complexity and practical relevance. More recently, Gonçalves et al. (2020) present a systematic literature review describing the history and trends regarding the safety stock determination from an operations research perspective. They also highlight that the number of contributions to MEIO has seen a significant increase from the year 2005 onwards, and they list many potential directions and trends for future research. There are two widely known approaches in MEIO to determine safety stock levels: the stochastic-service model (Clark & Scarf, 1960) and the guaranteed-service model (GSM), introduced by Simpson (1958). Detailed comparisons between them can be found in De Smet et al. (2019), Graves and Willems (2003), and Simchi-Levi and Zhao (2012).

The objective of this work is to develop a MEIO model based on the GSM approach that accounts for issues and characteristics arising in industrial practice in order to provide an improved representation to support strategic decision-making. Many authors (e.g. Inderfurth, 1993; Minner, 1998; Eruguz et al., 2014) have developed extensions to the GSM, but to the best of our knowledge, nobody has developed a model that can achieve optimum safety stocks on complex supply chains while integrating all the features typical of industrial environments presented in this work. We also propose strategies to obtain efficient solutions.

This paper addresses the problem of a multi-echelon, multi-product supply chain with both demand and lead time uncertainty. Demand can occur at any node in the network. This can result in hybrid nodes that have both dependent and independent demands. To the best of our knowledge, these characteristics, which represent the common operation mode of many multi-echelon systems, has not been addressed before, as most of the literature on supply chain inventory management considers only external demand at the final nodes of the network. Second, manufacturing plants can be placed at any location in the network, enabling the manufacture of any desired good at those locations. This feature allows generalizing and managing larger supply chains that have grown in their vertical integration. Capturing wider networks can significantly improve the inventory decision-making process across the process supply chains as is seen in those industries that produce both raw materials and finished goods. Third, the fill rate, which is the most widely applied service level measure in industry (Teunter et al., 2017), can be used as an alternative customer service indicator when setting safety stock levels. We adapt the fill rate constraint (Axsäter, 2006; Chopra & Meindl, 2013) to include hybrid nodes, and we propose a quadratic regression to estimate the equivalent Cycle Service Level (CSL) when fill rates required used in the model. In addition, Minimum Order Quantities (MOQ) for replenishment orders are explicitly modelled. Finally, the resulting nonconvex Nonlinear Programming (NLP) model is reformulated as a Quadratically Constrained Problem (QCP) by exploiting the structure of the constraints of the base model. Several computational examples for illustrative and industrial systems are presented to illustrate the application of the proposed model and its resulting improved computational performance.

The outline of the paper is as follows. The literature review and background with the basic concepts of the GSM are presented in the following subsection. The problem statement is given in Section 2, followed by the model formulation in Section 3. Section 4 details the application of the model on illustrative and real-world case studies. We conclude this article in Section 5. A Nomenclature section is presented at the end to facilitate the model understanding. A Supporting Information Section is included to provide the data input used in the real case study and detail additional discoveries relating to the impact of MOQ on service metrics.

## 1.1 Literature and Background of the Guaranteed-service Model

The present paper relies on the GSM to optimize safety stocks. Although this approach was developed more than 50 years ago, 80% of the existing works on this topic have been published in the last 2 decades (Eruguz et al., 2016). The first multi-echelon serial system for the GSM model was proposed by Simpson (1958), and then it was extended to deal with different network topologies (Graves & Willems, 2000; Inderfurth, 1993; Inderfurth & Minner, 1998; Minner, 1998). Later, Magnanti et al. (2006) developed a guaranteed-service approach for general acyclic networks. The main idea of the classic guaranteed service approach is that if the customer places an order of size $d_j(t)$ on node $j$ at time $t$, this order will be fulfilled by time $t + S_j$ (Graves & Willems, 2000), with $S_j$ being the guaranteed-service time of node $j$. Moreover, each node $j$ receives a service commitment from its upstream node $i \in J$, called inbound service time $SI_j$ ($SI_j = S_i$), and has an order processing time or lead time of $LT_j$. This lead time represents the time until the outputs are available to serve the demand, including material handling and transportation times. Both $SI_j$ and $LT_j$ are times that must be taken into account to define $S_j$. The Net Lead Time (NLT) is a concept that links them and represents the period of exposure that is not covered within the guaranteed service time and must be covered with safety stock. The *NLT* for node $j$ is defined as $NLT_j = SI_j + LT_j - S_j$. Figure 1 displays examples for different values of $S_j$. The first example (1) is the case where node $j$ promises to its customer a guaranteed service time equal to the worst-case replenishment time ($S_j = SI_j + LT_j$). This node places an order to its predecessor every time it receives an order from its customer, then it waits for the upstream node to process its order before processing the order without storing any inventory. In this case, $NLT_j = 0$. On the other hand, if the customer bounds the maximum possible service time ($S_j \leq maxS_j$) and this maximum is less than the worst-case replenishment time ($maxS_j \leq SI_j + LT_j$), node $j$ should satisfy customer demand in less time that the required to place an order on the supplier and process it. Therefore, $NLT_j > 0$, meaning that there is a period of time that should be covered with safety stock, as shown in cases (2) and (3) in Figure 1(A).

The objective function in the GSM is the total holding cost minimization. The holding cost on a given echelon $h_j$ is multiplied by the safety stock on that echelon. Assuming normal distribution to represent external demand patterns, the safety stock of echelon $j$ is calculated as $SS_j = k_j \, \sigma_j \sqrt{NLT_j}$, where $k_j$ is the safety stock factor that reflects the percentage of time that the safety stock covers the demand variation, and $\sigma_j$ represents the demand standard deviation. More details on the safety stock formula can be found in Eruguz et al. (2016). The aim of the GSM is to define the values of $SI_j$ and $S_j$ in order to reduce the safety stock holding cost. As shown in Figure 1(B), the guaranteed-service time defined for one node impacts the downstream stages in the network, because the guaranteed service time for the node becomes the inbound service time for its downstream successors ($SI_j = S_i$). In case (1), avoiding safety stocks in node $j$ yields large inventory levels on the successor stage $k$ (proportional to $NLT_k$), while in case (2) the inventory level at $k$ is reduced by holding stock in $j$.

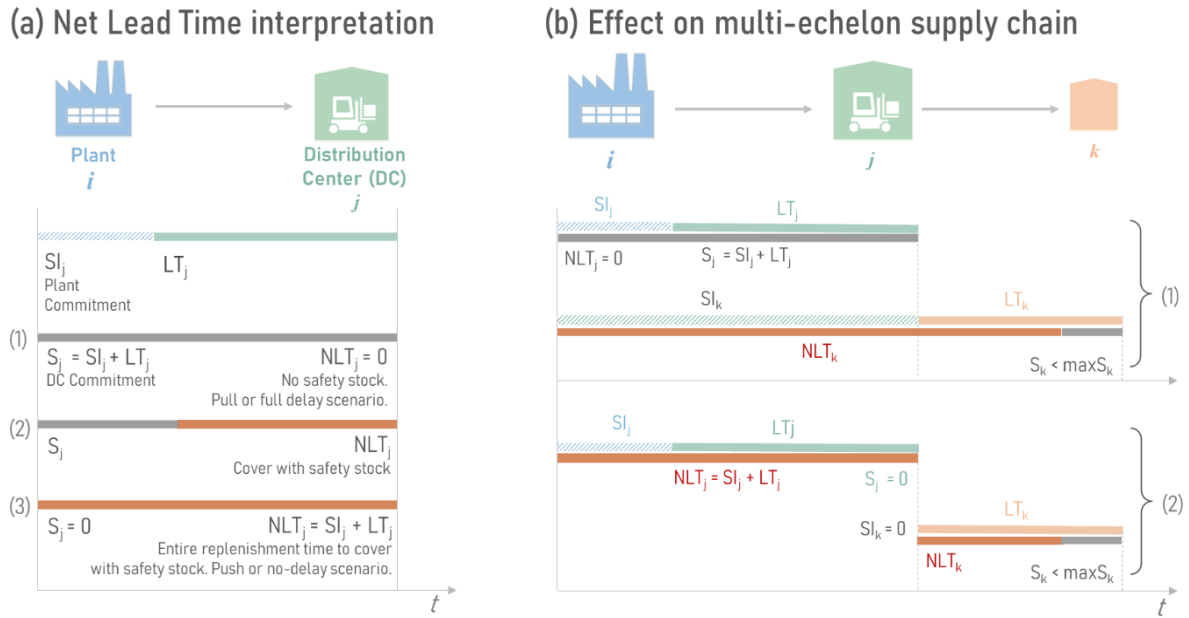

Figure 1: (a) Examples of different values for $SI_j$, $S_j$, and $NLT_j$: 1) pull scenario, 2) intermediate scenario, and 3) push scenario. (b) Guaranteed service approach in multi-echelon supply chains.

A key assumption of the basic GSM is that demand is bounded. If demand in a certain period exceeds the bound, it is assumed that other extraordinary measures such as overtime production are used to satisfy excess demand. Moreover, it is assumed that each stage of the supply chain operates under a (*R,S*) inventory policy with a base-stock level. The demand is independent and identically distributed following a normal distribution. Lead times are constant, and independent demand only occurs at final nodes in the network. In addition, the service times at the initial and the final nodes are inputs. Finally, each plant has a coefficient that represents the Bill of Materials (BOM) for product transformation and depends on location-location relationships (network arcs).

Over the years, many authors have worked on extending the original GSM assumptions to enable real-world supply chain characteristics to be captured, as presented in the survey by Eruguz et al. (2016). The authors in this survey summarize several extensions made to the basic model of Graves and Willems (2000). The main assumptions that were relaxed are related to the external demand, lead time variability, capacity constraints, service time customization, alternative replenishment policies, review periods, and extraordinary measures. Moreover, other authors have presented works about integrating the classic GSM with other activities or approaches. You and Grossmann (2008) develop models and algorithms that simultaneously consider inventory optimization and supply chain network design under demand uncertainty. In a subsequent work, these authors present an integrated multi-echelon supply chain design and inventory management model under uncertainty using the GSM (You & Grossmann, 2009). Klosterhalfen et al. (2013) propose an integrated hybrid guaranteed-service and stochastic-service approach for inventory optimization, that allows selecting the approach that minimizes costs. Recently, the work by Ghadimi (2020) presents a model for joint optimization of production capacity and safety stocks under the GSM approach. Bendadou et al. (2021) analyze the impact of merging activities in a supply chain under the GSM.

In addition to the extensions mentioned in (Eruguz et al., 2016), the inclusion of MOQ and fill rate as a service level measure have significant importance for representing supply chain dynamics and must be accounted for. Chopra and Meindl (2013) define the Cycle Service Level (CSL) as the fraction of replenishment cycles that end with all the customer demand being met, where the replenishment cycle is the interval between two successive replenishment deliveries. On the other hand, the product fill rate (*fr*) is the fraction of product demand that is fulfilled on time from the product in inventory. Chopra and Meindl (2013) describe how to introduce the fill rate given a continuous review inventory policy with a formula that links both indicators to obtain the equivalent CSL for single echelon networks. They also describe how a large MOQ yields larger fill rates. Silver and Bischak (2011) present an exact fill rate in a periodic review base stock system under normally distributed demand, and they state that the fill rate depends on four parameters, safety factor, coefficient of variation, review period, and lead time, but not on the minimum order quantity. De Smet et al. (2019) combine stochastic lead times with batching decisions for a distribution network based on the work of Humair et al. (2013) and calculate fill rates with an iterative procedure. More recently, Peeters (2020) accounts for MOQ to set safety stock levels and review periods integrated with stochastic lead times, based on the approach proposed by Humair et al. (2013), and using the Cycle Service Level (CSL) as a customer service measure.

In summary, this paper extends the GSM approach by including four main contributions. First, we address a more general supply chain that is frequently found in the pharmaceutical industry, with hybrid nodes including differentiated service times for each type of demand. We should note that this type of structure has not been considered previously in the literature, and it has great relevance in industrial practice for solving industrial problems. Second, we combined new (e.g. hybrid nodes) and existing features (e.g. stochastic lead times, review periods, fill rate) into a single model, requiring adaptations to allow their integration, such as extending the approaches of Inderfurth (1993) of stochastic lead times or Eruguz (2014) to a more generalized network. Third, we propose a new approximation method to include the fill rate as a target, using polynomial regression and adapting the existing formula in Chopra & Meindl (2013) for the GSM (*R,S*) inventory policy with minimum order quantities, multi-echelon networks, and hybrid nodes. A high R-squared value is obtained from the regression, meaning that the approximation has a good fit and is considered to be reliable. Finally, we introduce an exact reformulation of the NLP problem to a QCP to improve the computational efficiency by several orders of magnitude. This reformulation is equivalent to the original NLP problem, yielding the same optimal solution, thus guaranteeing the same quality of the solution. The proposed model is tested on examples for illustrative and industrial systems and provides computationally efficient solutions. The simulation of the results shows the accuracy of the proposed model to meet the service levels in the multi-echelon system under study.

## 2. PROBLEM STATEMENT

We are given a supply chain with a fixed design for a set of materials $p \in P$ that can be either raw materials or finished goods. The locations $j \in J$ belong to a set of plants, distribution centers, and retailers that can store different materials. Stock holding costs are incurred at all nodes; their unit costs are given. We assume uncertain demand and lead times. The objective is to determine the guaranteed-service times for each material at each location, and consequently how much safety stock to maintain at each location to minimize the total holding costs and satisfy a specified customer service level.

Unlike most literature on the topic, this work does not assume there is a final customer demand zone. In practice, it is usual that large hubs have an important external customer that places orders directly to this node. One the other hand, we assume that external or independent demand for any material $p$ can be placed at any node $j$ in the network. Each node $j$ can have an independent normally distributed demand of material $p$ with mean $\mu_{Ijp}$ and variance $\sigma_{Ijp}$, and/or an internal or dependent demand to satisfy replenishment orders from downstream nodes, with mean $\mu_{Djp}$ and variance $\sigma_{Djp}$. Demand is propagated upstream considering the risk pooling assumptions described in You and Grossmann (2009). A node that satisfies both dependent and independent demands is called a hybrid node, and an example of it is presented in Figure 2 on the left side.

Regarding the network topology, we assume divergent networks, as shown in Figure 2. In other words, a node that holds a material *p* can only receive this material from a single node and can distribute it to one or more locations, as is usual in finished goods supply chains. The same node can be supplied with another material $q \in P$ from another location, but this last one should be the only supplier of $q$ for that node. The route that each material follows, as well as the lead time distributions between two connected nodes, are given. Lead times are assumed to follow an independent normal distribution $LT_{jp} \sim N(lt_{jp}, \sigma_{LT_{jp}})$. They represent the delay that is under the responsibility of node $j$, including transportation, material handling, and other processing times until the material is ready to be shipped (i.e., is fulfilled).

Plants can be located at any node. Plant nodes can hold stock of both raw materials and finished goods. We introduce a general BOM based on a material-material relation, instead of a location-location relation as in Graves and Willems (2000). The value $\phi_{pq}$ determines the amount of material $p$ required to produce a unit of material $q$, regardless of the plant location. On the right side of Figure 2, there is an example of how independent, dependent and total demand mean and standard deviation of demand are propagated, including a Finished Good (FG) and a Raw Material (RM).

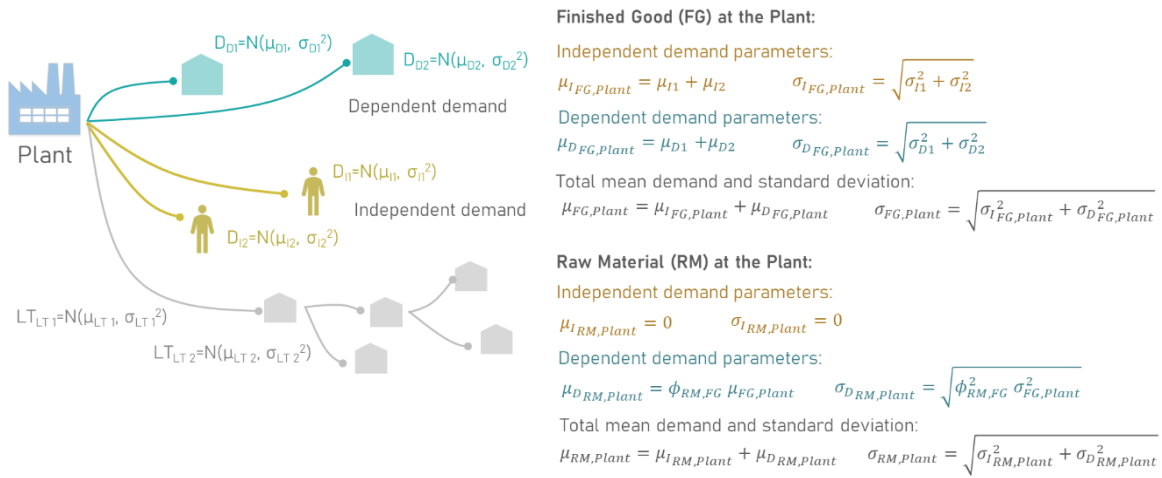

Figure 2: Example of a hybrid node and a divergent network.

We assume an $(R,S)$ inventory policy, with nested review periods ($r_{jp}$) as inputs to the model and common review days (Eruguz et al., 2014). Furthermore, a minimum order quantity $moq_{jp}$ may be enforced on replenishment orders. This means that if location $j$ needs to place an order, it will need to order at least the $moq_{jp}$, which may force it to receive an amount larger than required.

The network topology and the modes of transportation are assumed to be fixed. Hence transportation costs are not accounted for in this study. The service time of the most upstream nodes in the network and the maximum service time of each final node are given. We assume that information about demand is shared in the network and ordering decisions are decentralized; therefore, each node makes its own replenishment decisions and has no delay in ordering. For each node, the safety factor $k$ related to the CSL, which is represented by the standard normal distribution is also given, reflecting the percentage of time that the safety stock covers the demand variation. Alternatively, the modeler can also ask for a fill rate to be considered as a target service measure.

## 3. MODEL FORMULATION

The multi-echelon safety inventory optimization problem can be formulated as a nonlinear program (NLP) that deals with the safety inventory planning in a given supply chain. The model proposed in Graves and Willems (2000) is used as a basis, and all sets, parameters, and variables of this model are presented in the Nomenclature section. First, we assume that external demand on each node is a random variable defined as in Graves and Willems (2000). We assume that the external demand is normally distributed $D_{Ijp} \sim N(\mu_{I_{jp}}, \sigma_{I_{jp}})$, and its mean and standard deviations are propagated upstream to define internal demand means and standard deviations throughout the network. If there are stages with more than one successor, we require a decision on how to combine the demand bounds for the downstream stages to obtain a relevant demand bound for the upstream stage to position the safety stock (Graves & Willems, 2000). There will be a relative reduction in variability as we combine demand streams due to risk pooling. Therefore, the dependent demand parameters for material $p$ at node $j$ are obtained by converting the demand parameters for all materials $q$ that require $p$ as an input at all successor nodes $k$. The conversion is done via the Bill of Materials (BOM) $\phi_{pq}$ as a pre-processing step. The set $\Phi$ contains all valid material transformations (from material $p$ to material $q$), i.e. the raw material $p$ is required for obtaining the finished good $q$. $A$ is a set with indices $(i,j,p)$ indicating that there is a feasible route for material $p$ from node $i$ to node $j$. Note that $q = p$ and $i \neq j$ if it is a distribution link ($i$ to $j$) of the same product $p$, and $q \neq p$ and $i = j$ if node $j$ is a plant location that produces $p$ from $q$. We assume that the total demand of a given node $D_{jp}$ is the sum of the orders placed by immediate successors $D_{Djp}$ plus any external orders $D_{Ijp}$, i.e. $D_{jp} = D_{Ijp} + D_{Djp}$. As random variables are assumed to be independent from each other, the total demand is a linear combination of normal distributed variables, being also normally distributed $D_{jp} \sim N(\mu_{jp}, \sigma_{jp})$, as stated in Graves and Willems (2000) and You and Grossmann (2008). Similarly, the dependent demand $D_{Djp} \sim N(\mu_{D_{jp}}, \sigma_{D_{jp}})$ is also a linear combination of demands of successors. Therefore, the total mean demand is the sum of the mean demands as shown in Equation (1), and the total demand standard deviation is calculated as in Equation (2). In this work, we include the first term in both equations referring to independent demand mean and standard deviation that can be placed at any node. Note that the second term on both equations is equivalent to the dependent demand mean and deviation, that is, $\mu_{Djp}$ and $\sigma_{Djp}$. We assume pooling of both types of demand parameters for propagation purposes. For nodes where material $p$ is distributed, rather than transformed into $q$, $p = q$ and $\phi_{pq} = 1$. For manufacturing nodes where material $p$ is transformed into material $q$, $\phi_{pq}$ equals the amount of $p$ consumed per unit of $q$.

$$\mu_{jp} = \mu_{I_{jp}} + \sum_{(p,q)\in\Phi} \sum_{(j,k,p)\in A} \phi_{pq}\mu_{kq} \qquad \forall j \in J, p \in P_j \qquad (1)$$

$$\sigma_{jp} = \sqrt{\sigma_{I_{jp}}^2 + \sum_{(p,q)\in\Phi} \sum_{(j,k,q)\in A} \phi_{pq}^2 \sigma_{kq}^2} \qquad \forall j \in J, p \in P_j \qquad (2)$$

### 3.1 Constraints

The first set of constraints is related to bounding the guaranteed-service time variables. Equation (3) defines the first inbound service time for the starting (source) nodes in the network $J^0$, where $si^0$ is a given input. Equation (4) links the inbound guaranteed-service time $SI_{jp}$ and the guaranteed-service time of its upstream node $S_{iq}$. If there is a maximum accepted delay for any material on a node, the inequality in (5) is active. In addition, Equation (6) fix the maximum accepted service time $S_{Ejp}$ exclusively for external demand nodes. This service times does not impact the inbound service time of downstream nodes, because it related to safety stock dedicated only for external customers.

$$SI_{jp} = si_{jp}^0 \qquad \forall\, j \in J^0, p \in P_j \tag{3}$$

$$SI_{jp} \geq S_{iq} \qquad \forall\, (i,j,p) \in A, (q,p) \in \Phi, p \in P_j \tag{4}$$

$$S_{jp} \leq maxS_{jp} \qquad \forall\, j \in J, p \in P_j \tag{5}$$

$$S_{E\,jp} \leq maxSE_{jp} \qquad \forall\, j \in J_p^I, p \in P_j \tag{6}$$

#### 3.1.1. Manufacturing locations

In this work, a manufacturing site has the possibility of storing both raw materials and finished goods in the same node *j*. To the best of our knowledge, this representation of the manufacturing site has not been addressed before. Despite the notation demonstrating that only one node is involved, it is possible to represent the plant as two artificial nodes connected by an arc that represents the manufacturing time, as depicted in Figure 3. If it is required, for example, that the safety stock of raw materials is enough to satisfy production demand immediately from stock, a maximum possible service time *maxS$_{jq}$ = 0* can be required. On the other hand, safety stocks of finished goods could be constrained as *maxSS$_{jp}$ = 0* if no stock is allowed in the manufacturing node. For this special case, *A* represents an enabled production process to obtain product *p* at node *j*. There is a production lead-time that can be constant or normally distributed with parameters *lt$_{jp}$* and *σ$_{LTjp}$* to represent the manufacturing time. In case external and internal customers are able to place orders of finished goods to the plant, there would be dedicated safety stocks for each type of demand. The demand parameters for raw materials are defined by Eqs. (1) and (2), with their corresponding coefficient *ϕ$_{qp}$* related to the BOM. For example, if the raw material *q1* at plant *j* is used to produce both materials *p1* and *p2*, the mean demand for this raw material is $\mu_{j,q1} = \phi_{q1,p1}\,\mu_{j,p1} + \phi_{q1,p2}\,\mu_{j,p2}$ and the standard deviation is $\sigma_{j,q1} = \sqrt{\phi_{q1,p1}^2\,\sigma_{j,p1}^2 + \phi_{q1,p2}^2\,\sigma_{j,p2}^2}$.

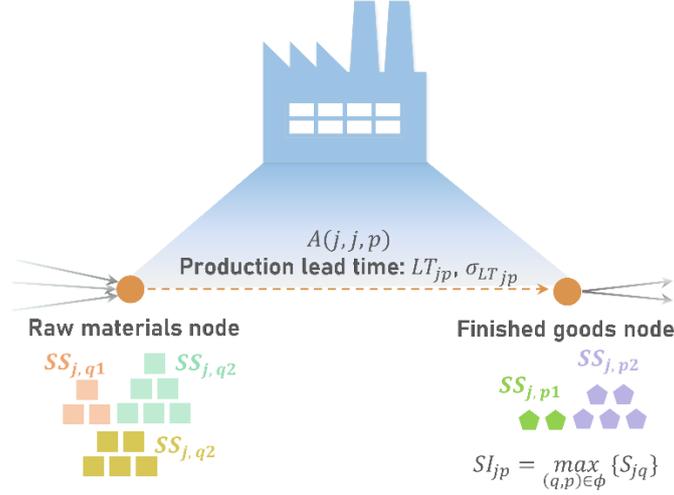

Figure 3: Representation of a manufacturing plant with optional safety stock for raw materials and finished goods.

### 3.1.2. *Stochastic lead times*

Concerning the incorporation of stochastic lead times to the GSM, our work is based on the approach by Inderfurth (1993). In that work, a serial network is proposed with external demand at the final nodes, called "demand nodes", and upstream nodes are called "non-demand nodes". Inderfurth (1993) proposes that the safety stock at a demand node involves the combination of the two random variables: independent demand random variable $D_{jp} \sim N(\mu_{I_{jp}}, \sigma_{I_{jp}})$ and the lead time random variable $LT_{jp} \sim N(lt_{jp}, \sigma_{LT_{jp}})$. The service process is usually arranged in such a way that customer demand will be fulfilled as soon as the fluctuating final-stage lead time will allow it. Therefore, in this planning situation there are flexible lead times. In our work, we propose to replace the mean lead time with the value of the net lead time variable. Hence, the safety stock for demand nodes can be calculated as:

$$SS_{I_{jp}} = k_{jp} \sqrt{NLT_{jp}\, \sigma_{I_{jp}}{}^2 + \mu_{I_{jp}}{}^2\, \sigma_{LT_{jp}}{}^2} \qquad \forall\, j \in J_p^I,\ p \in P_j \qquad (7)$$

where $SS_{I_{jp}}$ represents the safety stock level to satisfy independent demand. On the other hand, the approach of Inderfurth states that in an upstream stage the stochastic lead time is converted into a deterministic lead time to be consistent with integrated multi-level production planning from MRP systems. Therefore, $\widehat{lt}_{jp} = lt_{jp} + k_{LT_{jp}}\, \sigma_{LT_{jp}}$, where $k_{LT_{jp}}$ relates to the service level that denotes the probability regards the random lead time does not exceed the planned lead time $\widehat{lt}_{jp}$. Thus, using fixed planned lead times means that each internal demand will be satisfied just after $\widehat{lt}_{jp}$ periods predicting the expected stocks as follows:

$$SS_{D_{jp}} = k_{jp}\, \sigma_{D_{jp}} \sqrt{SI_{jp} + \widehat{lt}_{jp} - S_{jp}} \qquad \forall\, j \in J_p^D,\ p \in P_j \qquad (8)$$

Note that in the last case the pipeline inventory, $\mu_{jp}\widehat{lt}_{jp}$, differs from the deterministic case. In our work, we propose that the safety stock for each hybrid node is equal to the summation of safety stocks for independent and dependent demands, $SS_{jp} = SS_{I_{jp}} + SS_{D_{jp}}$. Hence, on each node, we define a safety stock to satisfy downstream orders and another safety stock for external orders. In general, these inventories are at the same location but have to be dedicated to each type of demand. Therefore, this limits pooling at this location. In practice, safety stock can be considered as a whole to satisfy demand if it does not mean a stockout on other customers at the same location.

### 3.1.3. Review periods

The original GSM assumes that review periods are common for all stages, and that the lead time includes any waiting and processing time at the given stage (Graves & Willems, 2000). Inderfurth (1993) assumes that a final-stage lead time additionally contains the length of the review period. In the present work, we introduce a more detailed definition of how review periods are accounted for in the net lead time definition. In Equations (7) and (8), no details are given on how the lead time mean is determined. In other words, these equations account for lead time variability, but it is not certain what $lt_{jp}$ includes in this delay time. We propose to distinguish this feature separating the review period $r_{jp}$ from the lead time parameter $lt_{jp}$. We assume nested review periods with common review days as in the work of Eruguz et al. (2014), and that a replenishment order is ready to satisfy the demand on its period of arrival. A node that faces external demand, needs to cover with safety stock the demand during net lead time $NLT_{jp} = SI_{jp} - S_{Ejp} + lt_{jp} + r_{jp}$. Hence, when an order is placed in a node, it is instantaneously propagated upstream. Therefore, in upstream nodes $NLT_{jp} = SI_{jp} - S_{jp} + lt_{jp} + r_{jp} - 1$, as stated in Eruguz et al. (2014).

As mentioned above, we propose the alternative of hybrid nodes with both types of demand. Inequalities (9) and (10) account for the definition of the net lead times to be covered with safety stock to achieve the desired service level for independent and dependent customers, respectively. These equations combine review periods and the stochastic lead time approach developed above. Note that $SI_{jp}$, $lt_{jp}$ and $r_{jp}$ are assumed to be the same for both types of demands, and $ARG_1$ and $ARG_2$ are positive continuous variables representing the terms in the square roots for the independent and dependent customers safety stocks, respectively. We can have different amounts of safety stocks for satisfying internal or external demands, because stochastic lead times are accounted for differently for both demand types.

$$ARG_{1\,jp} \geq SI_{jp} - S_{jp} + lt_{jp} + k_{LT_{jp}} \sigma_{LT_{jp}} + r_{jp} - 1 \qquad \forall j \in J^D, p \in P_j \qquad (9)$$

$$ARG_{2\,jp} \geq (SI_{jp} - S_{E_{jp}} + lt_{jp} + r_{jp}) \sigma_{I_{jp}}^2 + \mu_{I_{jp}}^2 \sigma_{LT_{jp}}^2 \qquad \forall j \in J^I, p \in P_j \qquad (10)$$

Note that the the right-hand sides in (9) and (10) must be positive. For this purpose, the upper bounds $S_{jp}$ and $S_{Ejp}$ are defined with inequalities (11) and (12). The former defines the upper bound as only for the case of dependent demand, while the latter one accounts for the upper bound in the case of independent demand.

$$S_{jp} \leq SI_{jp} + r_{jp} - 1 + lt_{jp} + k_{LT_{jp}} \, \sigma_{LT_{jp}} \qquad \forall j \in J^D, p \in P_j \qquad (11)$$

$$S_{E_{jp}} \leq SI_{jp} + r_{jp} + lt_{jp} + \left(\frac{\mu_{I_{jp}}}{\sigma_{I_{jp}}}\sigma_{LT_{jp}}\right)^2 \qquad \forall j \in J^I_p, p \in P_j \qquad (12)$$

### 3.1.4. Fill rate as a target service level

As described previously, the GSM uses the Cycle Service Level (CSL) as the customer service performance indicator when setting safety stocks. Since fill rate is more widely used in industry (Teunter et al., 2017), we extend the GSM to allow specifying fill rates if desired. Fill rates represent the fraction of demand that was met on-time from inventory. We use the works of Axsäter (2006) and Chopra and Meindl (2013) as a baseline and we propose modifications to account for additional features. First, we propose to replace the lead-time demand variability, expressed as $\sigma_{jp}\sqrt{NLT_{jp}}$, by $\sigma_{D_{jp}}\sqrt{ARG_{1_{jp}}} + \sqrt{ARG_{2_{jp}}}$, allowing the representation of multi-echelon networks with service times, stochastic lead times and hybrid nodes. This is possible since in Section 3.1.2 it is assumed that the safety stock is the sum of independent and dependent safety stocks for that node, $SS_{jp} = SS_{I_{jp}} + SS_{D_{jp}}$, with common service level $k_{jp}$. From Equations (7), (8), (9) and (10) we obtain $SS_{jp} = k_{jp}\left(\sigma_{D_{jp}}\sqrt{ARG_{1_{jp}}} + \sqrt{ARG_{2_{jp}}}\right)$".

Moreover, in (R,Q) inventory policies, Q refers to the replenishment quantity, however, for periodic review policies this amount is variable. We assume an average replenishment quantity, $Q_{jp} = \mu_{jp} \, r_{jp}$. From the formula presented by Chopra and Meindl (2013) and including the extensions mentioned, we can obtain the constraint that links fill rate ($fr_{jp}$) to the safety factor $k_{jp}$, and consequently to the CSL. The safety factor $k_{jp}$ becomes a continuous positive variable $KV_{jp}$ for those materials and locations that have fill rate levels active. The objective is to find the lowest CSL level that can meet a defined fill rate, given by Equation (13). $F_s(KV_{jp})$ and $f_s(KV_{jp})$ correspond to the cumulative and density standard normal distributions functions, respectively. Therefore, the constraint proposed in this model to find the minimum CSL to achieve the desired fill rate is given by Equation (13).

$$fr_{jp} \leq \frac{\left(\sigma_{D_{jp}}\sqrt{ARG_{1_{jp}}} + \sqrt{ARG_{2_{jp}}}\right)}{Q_{jp}}\Big(KV_{jp}\left[1 - F_s(KV_{jp})\right] - f_s(KV_{jp})\Big) + 1 \qquad (13)$$

$$\forall j \in J, p \in P_j, (j,p) \in F$$

### 3.1.5. Minimum Order Quantity (MOQ)

This requirement is frequently found in practice, however, to the best of our knowledge, there is little literature that relates MOQ to safety inventories. When an MOQ is required, flexibility is reduced, because the customer needs to either order many units or none. However, this does not necessarily mean that the risk and safety stocks are increased. Figure 4 depicts the effect that MOQ has on inventories. Plot (a) presents inventory evolution through time for a periodic review policy with a review frequency of one week and no lead time ($lt_k = 0$). In grey color, we can see the safety stock level ($SS_k$) set to cover a proportion of the demand excess during the net lead time. The basestock level B denotes the order-up-to level that must be accounted for when a replenishment order is placed, equivalent to parameter S in the (R,S) inventory policy. The order quantity (Q) is equal to the expected mean demand during a review period ($\mu_{jp}r_{jp}$). Each replenishment cycle, that is, the time between two consecutive replenishments deliveries, has a probability of non-stocking out of 1 - α. The safety stock level is set to cover demand variability during the net lead time (1-α)100% of the times, this being the probability determined by the $k_{jp}$ factor.

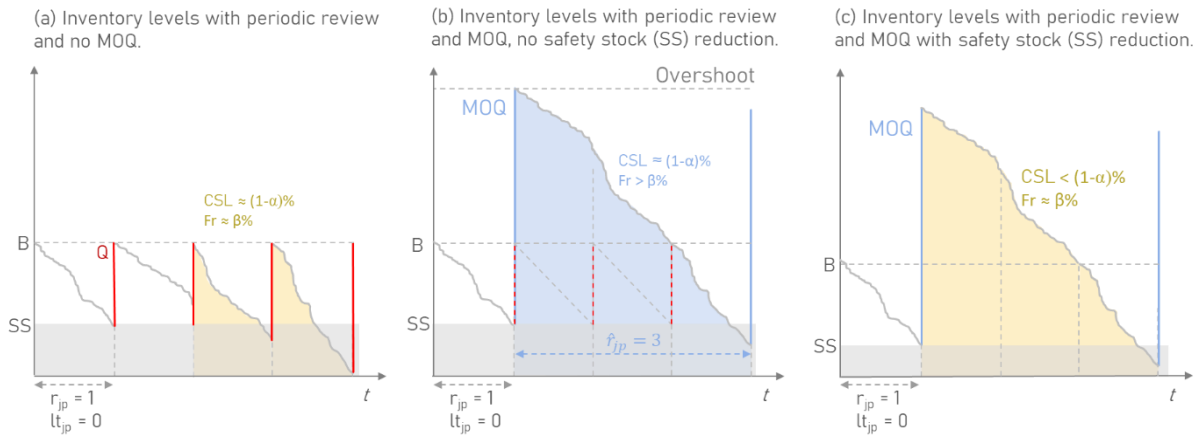

Figure 4: MOQ effect on inventory levels

If there is an MOQ required by the supplier to a node, and the MOQ is larger than the standard Q, the inventory level evolution will look like the one in Figure 4 (b). In this example, the MOQ size is three times the mean demand during the review period. Therefore, the first and the second periods have a low probability of stocking out, because there will be more inventory than is needed to satisfy the expected mean demand. However, the Cycle length is increased and a replenishment order is placed every three periods on average. The CSL measure will not be affected by the MOQ because there will be less replenishment cycles (Chopra & Meindl, 2013). On the other hand, fill rate levels will increase, which means that safety stocks can be reduced at the expense of increasing the cycle stock as a result of the MOQ requirement. The number of orders placed by the customer is not modified, and the overshoot in stock causes that many periods have more stock than necessary to fulfil the order. In this work, we propose to include this concept to reduce safety stock levels if the MOQ is larger than the original Q, as shown in Figure 4 (c). This reduction results in a decrease in the safety factor, because now $Q_{jp} = \max\{MOQ_{jp}, \mu_{jp}r_{jp}\}$ in (13). This feature represents and extension of the GSM to include an (s,S) inventory policy instead of the original (R,S) policy, in which s represents the Basestock level (B) and S the order-up-to level MOQ + SS. Figure 5 depicts how fill rates are generally larger than CSL for a given value of $KV_{jp}$. In the green lines, it is possible to see different curves for fill rates for increasing MOQ sizes, being *MOQ1* the smaller one and *MOQ4* the larger one. The larger the MOQ is, the larger is the fill rate achieved for a specified safety factor.

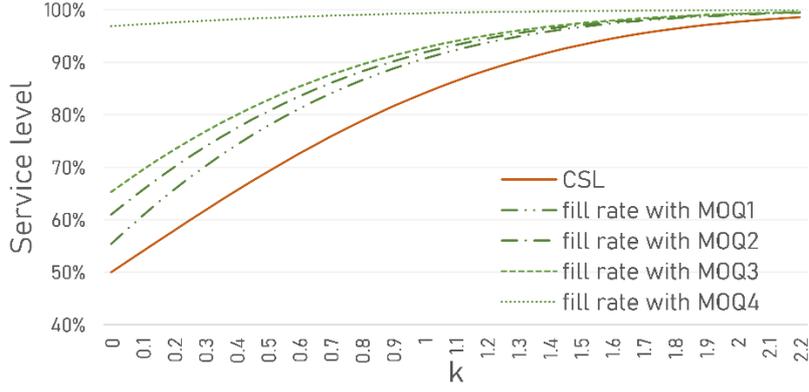

Figure 5: Fill rate sensitivity analysis with variations on replenishment quantities.

### 3.2 Objective Function

The objective function is to minimize safety stock holding cost as defined in Equation (14), where $h_{jp}$ is the coefficient that represents holding cost for each material $p$ at each location $j$.

$$min \sum_{j \in J} \sum_{p \in P_j} h_{jp} \, KV_{jp} \left( \sigma_{D_{jp}} \sqrt{ARG_{1\,jp}} + \sqrt{ARG_{2\,jp}} \right) \tag{14}$$

### 3.3 Solution Approach

The guaranteed service model (MNL1), given by equations (3)-(6), (9)-(14) is a nonconvex NLP with a concave objective function. Nonconvex NLP problems can in principle be solved with global optimization solvers like BARON. However, for medium or large-scale problem sizes, the computational time required to find the global optimum may be very expensive. To improve the model tractability and efficiency, we propose two solution approaches. In the first one, we use a quadratic regression to find an approximation to Equation (13), obtaining model MNL2. In the second one, we propose an exact reformulation of MNL2 into a quadratically constrained problem (QCP), denoted as MQC.

Equation (13) presents a difficulty to overcome. The function $g(x) = x\,[1- F_s(x)] - f_s(x)$ needs to be included in the mathematical model. To simplify this function, we propose a surrogate model through a second-order polynomial regression ($h(x) = ax^2 + bx + c$) to generate an approximation to $g(x)$ in (15), using as the domain the values that variable $KV_{jp}$ can take. The best-fit values obtained for the parameters in $h(x)$ are $a = -0.074700$, $b = 0.331986$, $c = -0.357195$, with $R^2 = 0.98$. A high R-squared value like the one obtained means that the approximation has a good fit with a small error. Figure 6 presents the original function $g(x)$ and the surrogate model function $h(x)$. The mathematical nonlinear program obtained (MNL2) is a nonlinear program composed of Equations (3)-(6), (9)-(12), (14)-(15), differentiated from MNL1 by the approximation on the fill rate constraint.

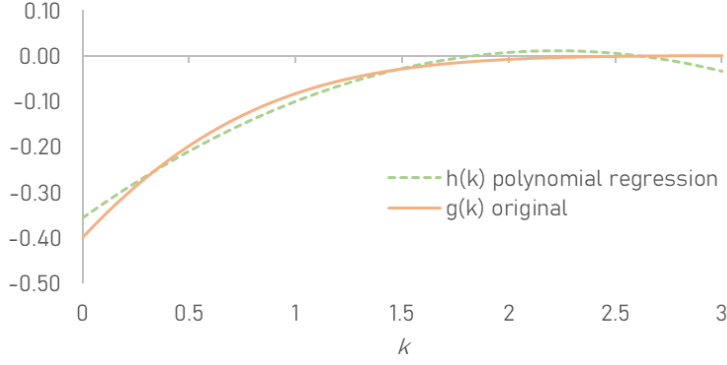

Figure 6: Surrogate model *h(x)* and original function *g(x)* curves.

$$fr_{jp} \leq \frac{1}{Q_{jp}} \left( \sigma_{D\,jp} \sqrt{ARG_{1\,jp}} + \sqrt{ARG_{2\,jp}} \right) \left( -a\, KV_{jp}^2 + b\, KV_{jp} - c \right) + 1 \qquad (15)$$

$$\forall\, j \in J,\, p \in P_j,\, (j,p) \in F$$

In order to improve the efficiency of the optimization, we propose a reformulation of the NLP model (MNL2) into a quadratically constrained problem, denoted MQC, which solvers like CPLEX and Gurobi can solve quite effectively in reasonable computational times. The idea behind the mathematical reformulation is to build an exact optimization model that benefits from its mathematical structure to improve the computational efficiency.

In order to derive the MQC reformulation, we first define a new variable $Z$ that replaces all the square root terms in the problem, where $Z = \sqrt{\tau}$ for a general expression $\tau$. Accordingly, the objective function of the NLP plotted in Figure 7(a) will be reformulated as in Figure 7(b), where Equation (16) is the reformulation of the objective function in (14).

$$\min \sum_{j \in J} \sum_{p \in P_j} h_{jp}\, KV_{jp} \left( \sigma_{D\,jp}\, Z1_{jp} + Z2_{jp} \right) \qquad (16)$$

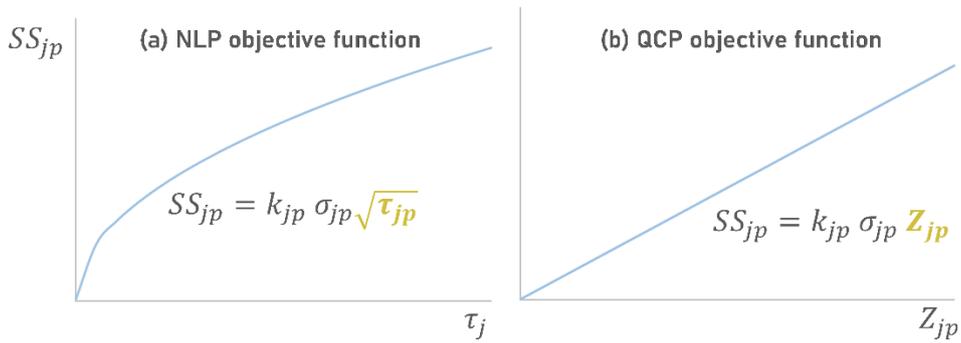

Figure 7: Objective function for NLP and QCP models.

Inequalities (9) and (10) are reformulated by replacing the left-hand sides terms with variables $Z1_{jp}$ and $Z2_{jp}$, resulting in Equations (17) and (18).

$$Z1_{jp}^2 \geq SI_{jp} - S_{jp} + lt_{jp} + k_{LT\,jp}\, \sigma_{LT\,jp} + r_{jp} - 1 \qquad \forall\, j \in J^D, p \in P_j \qquad (17)$$

$$Z2_{jp}^2 \geq (SI_{jp} - S_{E\,jp} + lt_{jp} + r_{jp})\,\sigma_{I_{jp}}{}^2 + \mu_{I_{jp}}{}^2\,\sigma_{LT_{jp}}{}^2 \qquad \forall\, j \in J^I, p \in P_j \qquad (18)$$

Finally, the fill rate constraint can be exactly reformulated as quadratically constrained, replacing Equation (15) with Equations (19) and (20).

$$Q_{jp}(fr_{jp} - 1) \leq \left(-a\,U_{jp} + b\,KV_{jp} - c\right)\left(\sigma_{D\,jp} Z1_{jp} + Z2_{jp}\right) \qquad (19)$$

$$\forall\, j \in J, p \in P_j, (j,p) \in F$$

$$KV_{jp}^2 - U_{jp} \leq 0 \qquad \forall\, j \in J, p \in P_j, (j,p) \in F \qquad (20)$$

In this way, the MQC reformulation is given by the objective function (16), subject to the constraints (3)-(6), (11)-(12), (17)-(20). This reformulation is equivalent to the original NLP as stated in the following proposition.

***Proposition 1.*** *The optimization problem MNL2 is equivalent to the optimization problem MQC.*

The proof of this proposition can be found in Appendix A.

## 4. APPLICATION AND RESULTS

### 4.1 Illustrative example and sensitivity analysis

An illustrative example is presented in Figure 8 to understand the model results and how different considerations impact safety stock decisions. From the supply chain network showed on the left, just a sample of products is selected, shown on the right. This case involves the production and distribution of a finished good (*SKU1*) obtained from two raw materials (*Raw1* and *Raw2*), and a Plant location that manufactures *SKU1* and delivers it to three retailers that satisfy external demand. The proportion of raw materials needed to obtain a unit of *SKU1* are $\phi_{Raw1,\,SKU1} = 1$ and $\phi_{Raw2,\,SKU1} = 0.014$. $S_{jp} = 0$ for *SKU1* at the three retailers' locations. The production lead time is 2 weeks and it is represented by the loop above the plant. Table 1 displays the demand and lead time input parameters, maximum service time constraints, and unit holding costs. The target CSL is 97% for all products.

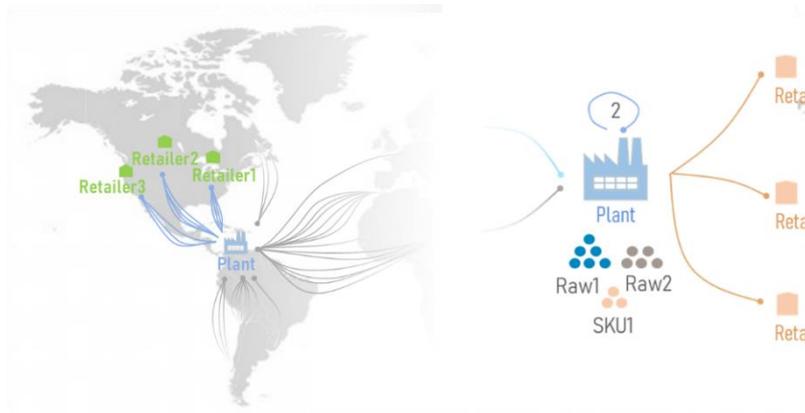

Figure 8: Illustrative example representation

Table 1: Illustrative example input data

| Material | | $Raw_1$ | $Raw_2$ | $SKU_1$ | $SKU_1$ | $SKU_1$ | $SKU_1$ |
|---|---|---|---|---|---|---|---|
| Location | | *Plant* | *Plant* | *Plant* | *Retailer$_1$* | *Retailer$_2$* | *Retailer$_3$* |
| Demand | $\mu_{jp}$ | 425,717 | 5,913 | 425,717 | 162,379 | 67,284 | 196,054 |
| (units) | $\sigma_{jp}$ | 192,229 | 2,669 | 192,229 | 119,665 | 61,585 | 137,258 |
| Coefficient of Variation (CV= $\sigma_{jp}/\mu_{jp}$) | | 0.45 | 0.45 | 0.45 | 0.73 | 0.91 | 0.70 |
| Lead Time | $LT_{jp}$ | 6 | 3 | 2 | 1 | 1 | 1 |
| (weeks) | $\sigma_{LTjp}$ | 1.9 | 0.7 | 0.0 | 0.3 | 0.6 | 0.4 |
| Max Service Time $S_{jp}$ (weeks) | | - | - | - | 0 | 0 | 0 |
| $h_{jp}$ ($/unit) | | 0.01171 | 0.00002 | 0.12 | 0.12 | 0.12 | 0.12 |

Results are detailed in Table 2. The computational tests are performed on an Intel® Core i7 CPU with 16 GB RAM and 4 parallel threads using Gurobi 9.1.2 as the QCP solver. The model (MQC) involves 30 continuous variables and 34 constraints. The CPU time required to obtain the optimal solution is 0.15 seconds and the total holding cost obtained is $162,205. It is possible to see that in the plant the decision is not to hold safety stock, and to select a guaranteed service of 2 weeks for supplying the retailers.

Table 2: Illustrative example results

| Material | $Raw_1$ | $Raw_2$ | $SKU_1$ | $SKU_1$ | $SKU_1$ | $SKU_1$ |
|---|---|---|---|---|---|---|
| Location | *Plant* | *Plant* | *Plant* | *Retailer$_1$* | *Retailer$_2$* | *Retailer$_3$* |
| $S_{jp}$ (weeks) | 0 | 0 | 2 | 0 | 0 | 0 |
| $SS_{jp}$ (units) | 1,143,300 | 11,228 | 0 | 459,359 | 243,783 | 536,961 |
| Holding cost$_{jp}$ ($) | 13,393 | 0.2 | 0 | 55,123 | 29,254 | 64,435 |

It is worth mentioning that the guaranteed service time of $SKU_1$ in *Plant* affects the retailers' safety stock levels, which need to cover for 2 more weeks with stock, as this is the inbound service time ($SI_{Retauler,SKU1}$ = 2 weeks). If this inbound service time continues to increase, the safety stock at the retailers will also increase, and it is possible that the model decides to change safety stock setting in the plant so as to take advantage of system-wide risk-pooling, and have a lower cost in the supply chain. As a sensitive analysis, Figure 9 (A) depicts the current case, with a production lead time of 2 weeks. In case (B), the production lead time is increased to 10 weeks, and the optimal solution for this case is to pool, holding stock of finished goods in the plant, with a total holding cost of $259,250. The decision of pooling in the plant yields a lower cost than if we maintained the decision of no safety stock in the plant

for *SKU1*, as shown in (C). The total cost increased to $265,360 in comparison to the solution in (B) because the opportunity of pooling and reducing the inbound service time for retailers is missed.

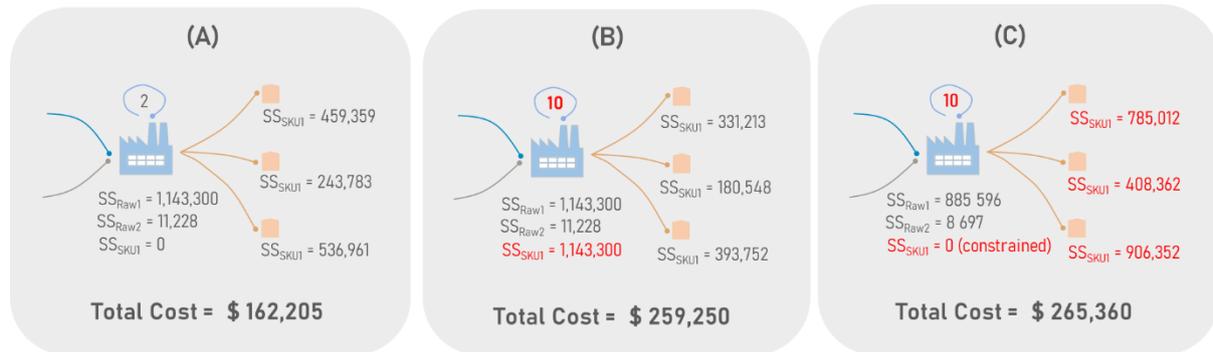

Figure 9: Analysis different lead times: A) 2-week production lead time, B) 10-week production lead time, C) 10-week lead time without safety stock for *SKU1* at the *Plant* as constraint.

This network does not have any hybrid node. As an example, we assume that the plant has external demand that is equal to the external demand of *Retailer₁*, with $\mu_{jp}$= 162,379 and $\sigma_{jp}$ =119,665. In that case, the optimal solution defines $SS_{plant,SKU1}$ = 459,359, and this stock will be set exclusively to satisfy external demand orders. Internal demand will still be satisfied using a MTO policy. Regarding raw materials, the safety stocks increase to face the increased demand the plant of the finished good, being $SS_{plant,raw2}$ = 13,632 and $SS_{plant,raw1}$ = 1,384,128.

For this illustrative example, there is no MOQ requirement. Therefore, we assume that for a periodic review policy, the expected order size ($Q_{jp}$) is equal to the mean demand during the review period, that is $\mu_{jp} r_{jp}$. The last analysis of this illustrative example concerns the measurement of customer service. In the current case, there is a desired 97% CSL, which corresponds to a safety factor *k* of 1.88. If a location, for example, *Retailer 1*, has to change the customer service metric from CSL to fill rate for a given product, a different CSL can be required to achieve the expected fill rate, so $k_{jp}$ (now $KV_{jp}$) is variable. In Figure 10, CSL and safety stocks are given for different cases varying target fill rates and MOQ constraints for *SKU1* on *Retailer 1*. The blue line indicates the expected fill rate, which is a given input in all cases except in the first one, in which the CSL is defined to set safety stocks as the original example, and the fill rate in this case is obtained through (19). In the following scenarios the target is a given fill rate (98%, 90%, 80% and 70%), and the CSL is obtained by the model. The yellow dashed line represents the resulting CSL for each case, and the yellow bar is the correspondent safety stock (secondary vertical axis) for that coverage. In addition, the brown dashed line and brown bars represent the resulting CSL and safety stock when applying an MOQ=500,000. The first case (the left-most case) is the current illustrative example scenario. The desired CSL is 97%, and a near 100% fill rate is expected, with or without a required MOQ. The second case sets a 98% fill rate to define safety stock levels. The minimum required CSL to achieve this expected fill rate decreases together with its corresponding safety stock level, and a sharper decrease occurs when a large MOQ is required. In the subsequent scenarios, the desired fill rates decrease and consequently, the CSL is lower. This difference is even more remarkable in the presence of MOQs, having no safety stocks defined for fill rates less than or equal to 80%, with a minimum required CSL of 50%.

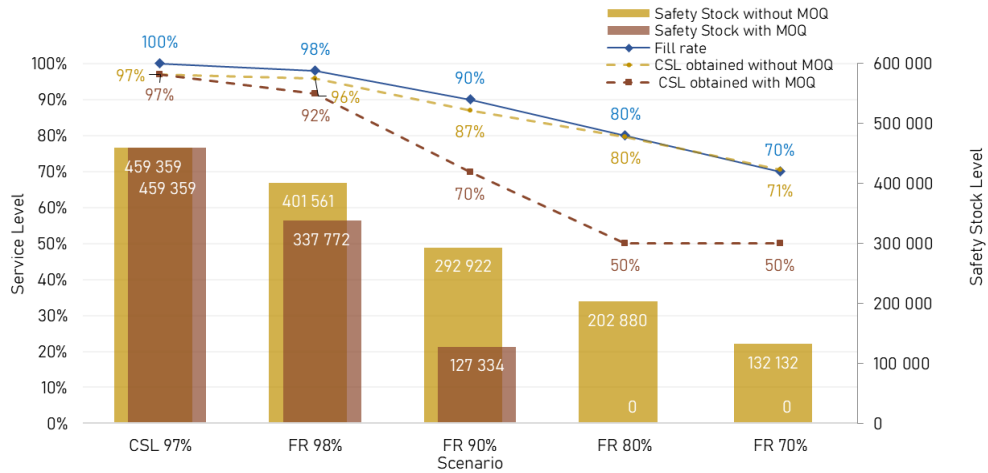

Figure 10: Effect of fill rate and MOQ on CSL and safety stock levels.

In summary, the illustrative example provides a clear demonstration about how the model is able to manage safety stock decisions in a multi-echelon network. The risk pooling allows the business to recognize the opportunity of potential savings by holding safety stock upstream given a lower total variability in the demand. Moreover, it is interesting to see the impact of a large minimum order quantity on safety stock. The larger the lot size, the less need for safety stock. However, it is worth to mention that it becomes an additional carrying cost for cycle stock. The MOQs represent transportation and production constraints that are frequent for found in almost all echelons in the supply chain. The opportunity of including this feature combined with the most used service level metric, the fill rate, yields significant savings in safety stock levels, even eliminating them in the case of large MOQs, as shown in Figure 10. In this figure it is also possible to see how the different service level metrics yield significantly different solutions, reinforcing the importance of representing the desired service level measure by the company.

### 4.2 Small-size industrial case study

The MQC formulation is now applied to a small-size industrial case study with the supply chain network shown on the left in Figure 8, with two echelons, 4 SKUs, and 31 raw materials coming from different locations and 1 intermediate product produced and consumed in the plant. The first echelon has one plant and the second echelon has three retailers, as in the illustrative example presented above. The complete input data is presented in Tables S1 and S2 in the Supporting Information Section. Note that lead times have decimals because they are averages of historic data. For MEIO purposes, the ceiling of the lead time $\lceil lt_{jp} \rceil$ is used as input.

The MQC model has 248 continuous variables and 291 constraints. While the NLP formulation using BARON is not able to find a feasible solution within 1000 seconds, using Gurobi the proposed

QCP formulation finds the optimal solution in 0.03 seconds. The optimal solution is $762,503. The results were compared with a commercial software (not identified due to confidentiality reasons), and the current safety stock levels for raw materials (RM) and finished goods (FG) are summarized in Table 3. While the commercial software solution obtains a 10% reduction in holding costs regarding the current safety stock levels, the proposed model in this work yields a 17% reduction, clearly showing the advantage of this tool to reduce the capital in inventory. Note that safety stocks at the retailers are slightly larger with the proposed model. The model seems to reduce the amount of inventory of raw materials, yielding the largest reduction in holding costs.

Table 3: Small-size industrial safety stock levels and holding costs

|  | Material | Model output | Baseline (Current level) | Commercial software |
|---|---|---|---|---|
| Holding cost | RM | $ 126,532 | $ 155,070 | $ 193,080 |
|  | FG | $ 635,971 | $ 761,525 | $ 634,950 |
| Safety stock level | RM | 8,031,692 | 9,917,801 | 13,029,797 |
|  | FG | 2,740,254 | 3,042 031 | 2,685,615 |

### 4.3 Medium-size industrial case study

This case involves the same network as in Figure 8, but now it has 20 finished goods and 120 raw materials, requiring 196 safety stock decisions. The MOQs constraints (Equations (19) and (20)) are active for all nodes and materials, and the customer service measure of interest is the fill rate, which is different for each material. The size of the QCP model is 1,973 constraints and 1,427 continuous variables, and is solved to optimality within 3 seconds using Gurobi as the QCP solver. This further supports the usefulness of the proposed approach for solving real-world problems efficiently, obtaining optimal solutions at low computational expense.

In summary, the results of the three case studies show the computational efficiency of the proposed model to obtain fast and optimal solutions for different instances. A significant reduction in computational time is of great relevance since the company is presently facing problems that can take a few days to run for 52 periods with commercial software. While the NLP formulation is not able to find a feasible solution, the proposed QCP formulation finds the optimal solution at minimum computational expense. This clearly shows the competitive advantage of this tool, while achieving customer service levels with a minimum capital in inventory. The simulations in the following section are helpful to test the accuracy of the model outputs to achieve expected service levels.

### 5. VALIDATION THROUGH SIMULATION

Having presented the model formulation and its application to several case studies, we include simulation studies to provide some insights about the effectiveness and the robustness of the predicted solutions.

## 5.1 Single-echelon simulations

Safety stock formulas are developed for normally distributed demands during lead time. Soares (2013) states that if the coefficient of variation (CV) is not considerably less than 1, there is a relatively high probability for negative demand when using the normal distribution, and the accuracy may be compromised. Generally, the cases addressed in this paper concern demand patterns that can be considered smooth in most of the cases, with relatively low CVs. The CVs in the cases addressed in this paper ranges from 0.45 to 0.91, while lead times CVs are within 0 and 0.6. For some values, such as 0.91, the service level target may be compromised, showing a slightly lower performance in reaching the expected service level. We propose to add a brief analysis to evaluate if the safety stock decisions are accurate to meet the service levels for different demand CVs and deterministic lead times. The complete analysis and results are presented in the Supporting Information, Section S2. Figure 11 displays the results on each plot for each service type. The average expected service level is on the horizontal axis and the one obtained from simulation is on the vertical axis. The orange line refers to the mean effective service level obtained from simulation, the light orange area represents the confidence intervals, and the grey line the ideal results. All 95% confidence intervals presented in Tables S7 and S8 in the Supporting Information Section for both service types include the target value within their bounds. Note that the estimation of CSL is slightly less accurate when the CVs are closer to 1, with a maximum difference of 0.02, and the estimation is better for low CVs, in general for CVs lower than 0.66. In the case of fill rates, only a few values differ in 0.01 or 0.02 points from the expected values, for CV 0.96. For low targets, the fill rate is larger than expected. In summary, with these simulations it is possible to see that, in general, the service levels are achieved with the safety stock levels proposed by the model, with a maximum difference on 0.02.

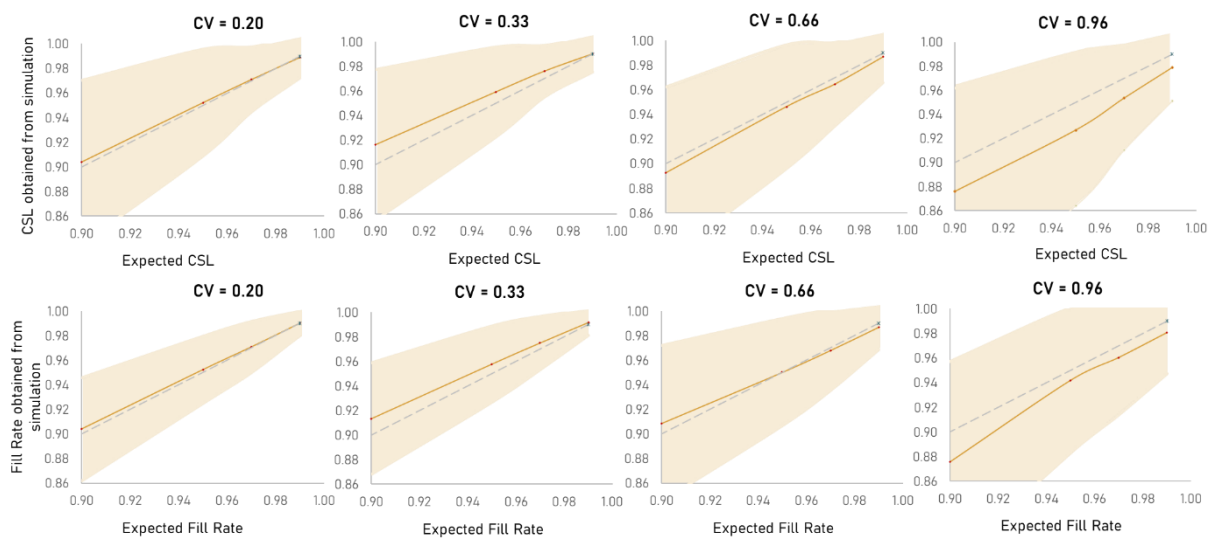

Figure 11: Expected vs. Effective service level obtained for different targets and CVs.

### 5.2 Multi-echelon simulations

The following simulations have the purpose to evaluate service level achievement in multi-echelon networks. All results obtained from the developed model are validated using simulation. This is done using an open-sourced discrete-time inventory simulation software package written in the Julia language: *InventoryManagement.jl* (Perez, 2021). This simulator allows modeling multiproduct supply networks of any topology (e.g., serial, divergent, convergent, tree, or general). Each of the features included in the extended GSM model can be simulated using this software: hybrid nodes, MOQ, bill of materials, stochastic demand, and stochastic lead times. For greater clarity, the validation of the illustrative example is presented with two extra scenarios to analyze how some features affect the system behavior in the simulation. The simulation tool also has a procedure to estimate the initial required parameters of normally distributed random variables, to obtain the desired parameters of the normal distribution after truncation of negative values.

Demand and lead times values are randomly generated using normal distributions for each period, with the parameters defined in Table 1. Random demand is only generated when there is external demand, and it is then propagated upstream. Base stock levels are calculated following Equation (21).

$$B_{jp} = SS_{jp} + \mu_{D_{jp}}\left(SI_{jp} - S_{jp} + \left\lceil lt_{jp} + z\,\sigma_{LT_{jp}}\right\rceil + r_{jp} - 1\right) + \mu_{I_{jp}}\left(SI_{jp} - S_{E_{jp}} + \left\lceil lt_{jp}\right\rceil + r_{jp}\right) \quad (21)$$
$$\forall\, j \in J_p,\, p \in P_j$$

In each period, a random demand value is generated, orders are placed and delivered, and inventory levels are updated. During each review period, an order is placed if the inventory position of a material on a location is below the basestock level. The management policy is decentralized: each location asks the amount they need to reach the basestock level, no matter how much the upstream node has on stock. If the available inventory is not able to meet demand, backorders are considered (the extraordinary measures that are referred to in the GSM approach are ignored in the simulation). The period selected in this case is one day with 7,000 days (1,000 weeks) in each run, so as to simulate the stationary state at each location, and 8 replications for each scenario are carried out. Demand and review periods have a weekly basis. The sequence of steps in the simulation of each period is the following:

1. External demand is observed and discounted at each node with independent demand. Unfulfilled demands are marked as lost sales.
2. On each review period, internal replenishment orders are placed if the current stock level is below the base stock level, and the lead time $lt_{jp}$ is defined accordingly for this order. Orders start being processed with no delay. Internal demand is discounted at each node, and unfulfilled replenishment orders are registered.
3. Stocks are updated with the replenishment orders that arrive at each node.

The illustrative case study was run and two additional scenarios were also tested for sensitivity analysis, combining fill rates as a target measure and MOQs for finished goods at retailers. Table 4 presents the model output for each scenario that is to be used by the simulation.

Table 4: Input data for simulation of all scenarios (multi-echelon)

| Scenario | Location | Material | Safety Stock level | Basestock level | Expected fill rate | Expected CSL | k factor |
|---|---|---|---|---|---|---|---|
| 1 | Plant | $SKU_1$ | 0 | 5,375,266 | - | 97% | 1.88 |
| | Plant | $raw_1$ | 1,118,096 | 39,992 | - | 97% | 1.88 |
| | Plant | $raw_2$ | 10,428 | 5,375,266 | - | 97% | 1.88 |
| | Retailer1 | $SKU_1$ | 459,166 | 1,108,682 | - | 97% | 1.88 |
| | Retailer2 | $SKU_1$ | 243,680 | 512,816 | - | 97% | 1.88 |
| | Retailer3 | $SKU_1$ | 536,736 | 1,320,952 | - | 97% | 1.88 |
| 2 | Plant | $SKU_1$ | 0 | 0 | 97% | - | - |
| | Plant | $raw_1$ | 910,989 | 5,168,159 | 97% | - | 1.58 |
| | Plant | $raw_2$ | 8,069 | 37,633 | 97% | - | 1.50 |
| | Retailer1 | $SKU_1$ | 383,857 | 1,033,373 | 97% | - | 1.57 |
| | Retailer2 | $SKU_1$ | 209,762 | 478,898 | 97% | - | 1.62 |
| | Retailer3 | $SKU_1$ | 446,787 | 1,231,003 | 97% | - | 1.56 |
| 3 | Plant | $SKU_1$ | 0 | 0 | 97% | - | - |
| | Plant | $raw_1$ | 910,989 | 5,168,159 | 97% | - | 1.58 |
| | Plant | $raw_2$ | 8,069 | 37,633 | 97% | - | 1.50 |
| | Retailer1 | $SKU_1$ | 301,155 | 950,671 | 97% | - | 1.23 |
| | Retailer2 | $SKU_1$ | 118,761 | 387,897 | 97% | - | 0.92 |
| | Retailer3 | $SKU_1$ | 369,736 | 1,153,952 | 97% | - | 1.30 |

The results of the scenarios are shown in Table 5. Results show the average of the effective service levels of each echelon obtained from the simulation and the 95% confidence intervals for these service levels.

The first scenario simulates the results of the illustrative example presented in Section 4.1, with a 97% CSL target for every material/location combination and no MOQs required. Note that service levels at the plant are larger than expected. As mentioned in (Minner, 1998) the approach of Inderfurth may induce large safety stocks, and does not benefit from joint coverage against both sources of uncertainty. On the other hand, as in the single echelon simulations, the effective CSL obtained in the retailers is slightly lower than the expected 97%. A possible reason is that the simulation model takes the ceiling of lead time values, what causes an increase in the lead time parameter. This error can be reduced if the discretization of time in the simulation is increased to simulate, for example, lead times given in hours. However, this impacts on the computational efficiency.

In Scenario 2, the service level target is changed to a 97% expected fill rate for all materials. Table 4 presents the base stock level and the minimum safety factor required to meet the target, obtained from variable $KV_{jp}$ in Equations (19) and (20) in the optimization model. Note that all safety factors are strongly reduced regarding Scenario 1, showing that lower safety stock levels are enough to meet the desired fill rate. The retailers can achieve in average the 97% target, and the service levels on upstream nodes are again larger than the target. The proposed quadratic regression to approximate the fill rate is accurate for the proposed scenario, according to simulation results.

Finally, Scenario 3 adds a minimum order quantity of 500,000 units to Scenario 2 to supply retailers. This minimum order quantity is the minimum necessary batch size to deliver an order from the plant to a retailer. The safety factor is reduced in Table 4 for the retailers because of the MOQ, as discussed in Section 3.1.5. Note that fill rates are achieved with a 24% of reduction in safety stock levels. The major reduction (43%) is in *Retailer 2*, the one that has the greatest CVs in demand and lead time. The MOQ effect can be detected on inventory levels through time in Figure 12, that presents inventory levels for *Retailer 2*. This retailer has the lowest mean demand, and therefore the MOQ is able to cover more review periods, 7 weeks on average. In the zoomed-in rectangle, it is possible to detect that an MOQ can cover several weeks of demand.

Table 5: Simulation Results

|          |          | Location     | Plant         |            | Plant         |            | Retailers     |            |
|          |          | Materials    | Raw materials |            | $SKU_1$       |            | $SKU_1$       |            |
| Scenario | Target   |              | Mean          | 95% CI     | Mean          | 95% CI     | Mean          | 95% CI     |
|----------|----------|--------------|---------------|------------|---------------|------------|---------------|------------|
| 1        | 97% CSL  | Effective CSL| 100.0%        | [100,100]  | 100.0%        | [100,100]  | 96.3%         | [95.9,96.8]|
| 2        | 97% FR   | Effective FR | 100.0%        | [99.9,100] | 100.0%        | [99.9,100] | 97.0%         | [96.4,97.6]|
| 3        | 97% FR   | Effective FR | 99.8%         | [99.4,100] | 99.8%         | [99.4,100] | 97.9%         | [96.9,98.9]|

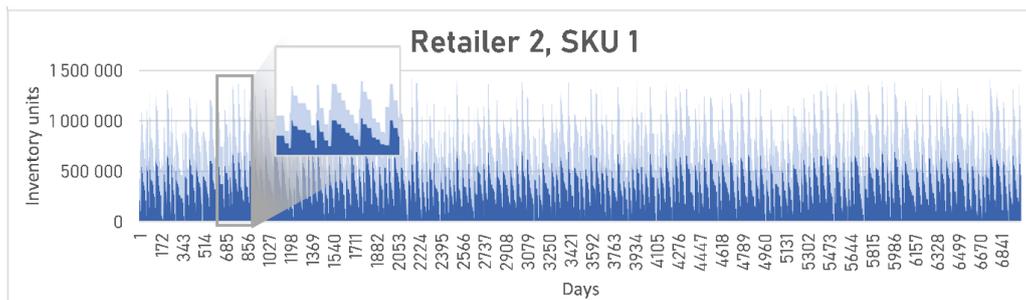

Figure 12: Inventory levels and inventory positions of *Retailer 2* from the simulation of Scenario 3

In summary, confidence intervals obtained from the simulation show the accuracy of the proposed model to meet the service levels in the multi-echelon system under study. The estimation is more accurate for fill rate targets. However, the differences between the expected and the effective CSLs is at most 0.02 points. The plant with a MTO policy for the finished goods, pushing more safety stocks to the retailers works well. The safety stocks of raw materials may be reduced, because they allow achieving larger service levels than expected.

## 6. CONCLUSIONS

In this paper, we have presented an optimization model based on the guaranteed-service approach that determines the optimal safety stock allocation in multi-echelon divergent networks. To the best of our knowledge, this is the first model that brings together multiple features typical of industrial practice, such as MOQs, hybrid nodes, and alternative service level measures to determine safety stock levels. It is also the first model to introduce the QCP reformulation to improve the computational efficiency of the optimization. The QCP outperforms the NLP formulation by allowing the use of QCP solvers, which leads to order of magnitude reductions in computational time. Real-world examples from the pharmaceutical industry are presented to illustrate the applicability of the proposed formulation. Optimal solutions can be found with small computational expense for medium/large scale problems. The simulation of the results demonstrates that the model is valid for achieving target service levels.

Results were presented through several exchanges with the pharmaceutical industry who had access to a commercial software vendor. The proposed model provides solutions with increased efficiency, apart from obtaining the exact global solutions. According to the feedback by the company, the results are significantly better than the current safety stock levels they currently have. We also obtain better performance than the commercial software used, which has missing features, such as hybrid nodes. A reduction in computational time is of great importance to the company, since they are presently facing solving problems that take days to run with the commercial software. In addition, they found valuable that the work provides an opportunity for developing the algorithm as an open source, and not hidden under a software package. In addition, it presents an opportunity to integrate MEIO decisions with other tactical planning model, such as rhythm wheel (lot size optimizer along with sequencing), which gives an end-to-end view of the optimization.

Future work will address an extension of the present formulation for cases of non-normal demand, and a pre-processing procedure of input data in order to decide which mathematical formulation is appropriate to optimally determine safety stock levels. The effects of CV and MOQ on CSL estimation can also be analyzed to review other potential safety stocks reductions. This research can also be extended by including responsive characteristics to account for supply chain disruptions and by including storage capacity limitations. Another important extension is constrained capacity on the nodes.

## 7. NOMENCLATURE

### 7.1 Sets

| | |
|---|---|
| $J$ | Set of locations |
| $P_j$ | Subsets of products that can be stored at location $j$ |
| $J_p$ | Subsets of locations in the route of material $p$ |
| $J^0$ | Subset of starting locations in the network |
| $J^I$ | Subset of locations that face external demand |
| $J^D$ | Subset of locations that face internal demand |
| $A$ | Subset of routes segments (from node $i$ to node $j$) enabled for material $p$ |
| $F$ | Set of locations that have materials with an active fill rate as a target |
| $\Phi$ | Set of all valid material transformations (from material p to material q) |

## 7.2 Parameters

| | |
|---|---|
| $\mu_{jp}$ | Mean of the total demand of material $p$ in location $j$ |
| $\sigma_{jp}$ | Standard deviation of the total demand of material $p$ in location $j$ |
| $\mu_{D\,jp}$ | Mean of the dependent demand of material $p$ in location $j$ |
| $\sigma_{D\,jp}$ | Standard deviation of the dependent demand of material $p$ in location $j$ |
| $\mu_{I\,jp}$ | Mean of the independent demand of material $p$ in location $j$ |
| $\sigma_{I\,jp}$ | Standard deviation of the independent demand of material $p$ in location $j$ |
| $lt_{jp}$ | Lead time/order processing time of material $p$ in location $j$ |
| $\sigma_{LT\,jp}$ | Standard deviation of the lead time/order processing time of material $p$ in location $j$ |
| $h_{jp}$ | Holding cost of material $p$ in location $j$ |
| $si^0_{jp}$ | Inbound service time for the source nodes in the network |
| $\phi_{pq}$ | Amount of material $p$ required to produce material a unit of material $q$ |
| $maxS_{jp}$ | Maximum service time accepted for material $p$ in location $j$ |
| $maxS_{E\,jp}$ | Maximum service time for material $p$ in location $j$ regarding external demand |
| $r_{jp}$ | Stock review period for material $p$ in location $j$ |
| $moq_{jp}$ | Minimum Order Quantity of material $p$ that location $j$ must place |
| $Q_{jp}$ | Replenishment order size of material $p$ at location $j$ |
| $fr_{jp}$ | Fill rate level of material $p$ at location $j$ |
| $k_{jp}$ | Safety factor associated with CSL of material $p$ at location $j$ |

## 7.3 Positive Variables

| | |
|---|---|
| $S_{jp}$ | Guaranteed service time within which location $j$ will attend demand of material $p$ |
| $S_{E\,jp}$ | Guaranteed service time for external demand of product $p$ at location $j$ |
| $SI_{jp}$ | Inbound Guaranteed service time at location $j$ of material $p$ |
| $ARG_{1\,jp}$ | Argument of square root for independent demand of material $p$ at node $j$ |
| $ARG_{2\,jp}$ | Argument of square root for dependent demand of material $p$ at node $j$ |

$NLT_{jp}$    Net Lead time of material $p$ at node $j$

$Z1_{jp}$    Variable used for quadratic reformulation on dependent demand net lead time formula

$Z2_{jp}$    Variable used for quadratic reformulation on independent demand net lead time formula

$KV_{jp}$    Variable used to replace $k$ input factor when the fill rate is introduced to determine safety stocks

$U_{jp}$    Variable defined to replace $KV_{jp}^2$ and avoid trilinear terms

## 8. APPENDIX A

**Proof of Proposition 1.** We define two positive continuous variables $Z1_{jp}$ and $Z2_{jp}$ as follows,

$$Z1_{jp} = \sqrt{ARG_{1\,jp}} \qquad \forall j \in J^D, p \in P_j \qquad (A.1)$$

$$Z2_{jp} = \sqrt{ARG_{2\,jp}} \qquad \forall j \in J^I, p \in P_j \qquad (A.2)$$

Substituting (A.1) and (A.2) into the objective function (Eq. (14)), we have:

$$\min \sum_{j \in J} \sum_{p \in P_j} h_{jp} \, KV_{jp} \left( \sigma_{D\,jp} Z1_{jp} + Z2_{jp} \right) \qquad (A.3)$$

Since (A.3) is a minimization problem and both variables are present in the objective function, following the KKT optimality conditions, we can relax Eqs. (A.1) and (A.2) and rewrite them as inequalities (A.4) and (A.5), both being active at the optimal solution.

$$Z1_{jp} \geq \sqrt{ARG_{1\,jp}} \qquad \forall j \in J^D, p \in P_j \qquad (A.4)$$

$$Z2_{jp} \geq \sqrt{ARG_{2\,jp}} \qquad \forall j \in J^I, p \in P_j \qquad (A.5)$$

(A.4) and (A.5) can be reformulated as quadratic inequalities,

$$Z1_{jp}^2 \geq ARG_{1\,jp} \qquad \forall j \in J^D, p \in P_j \qquad (A.6)$$

$$Z2_{jp}^2 \geq ARG_{2\,jp} \qquad \forall j \in J^I, p \in P_j \qquad (A.7)$$

Therefore, Eqs. (9), (10), and (15) can be rewritten as follows:

$$Z1_{jp}^2 \geq SI_{jp} - S_{jp} + lt_{jp} + k_{LT\,jp}\,\sigma_{LT\,jp} + r_{jp} - 1 \qquad \forall j \in J^D, p \in P_j \qquad (A.8)$$

$$Z2_{jp}^2 \geq (SI_{jp} - S_{E\,jp} + lt_{jp} + r_{jp})\,\sigma_{I_{jp}}{}^2 + \mu_{I_{jp}}{}^2\,\sigma_{LT\,jp}{}^2 \qquad \forall j \in J^I, p \in P_j \qquad (A.9)$$

$$fr_{jp} \leq \frac{1}{Q_{jp}}\left( \sigma_{D\,jp} Z1_{jp} + Z2_{jp} \right)\left(-aKV_{jp}^2 + bKV_{jp} - c\right) + 1 \qquad (A.10)$$

$$\forall j \in J, p \in P_j, (j,p) \in F$$

Moreover, we define $U_{jp} = KV_{jp}^2$ to avoid a trilinear term in (A.10). Following the same steps, we have variable $KV_{jp}$ in the objective function. Since (A.3) is a minimization problem, we rewrite the equality as the inequality:

$$U_{jp} \geq KV_{jp}^2 \qquad \forall j \in J, p \in P_j, (j,p) \in F \qquad (A.11)$$

Hence, $KV_{jp}^2 - U_{jp} \leq 0$, as in Eq. (20). Replacing $KV_{jp}^2$ with $U_{jp}$, we obtain Equation (A.12):

$$fr_{jp} \leq \frac{1}{Q_{jp}} \left( \sigma_{D_{jp}} Z1_{jp} + Z2_{jp} \right) \left( -aU_{jp} + bKV_{jp} - c \right) + 1 \qquad (A.12)$$

$$\forall j \in J, p \in P_j, (j,p) \in F$$

Therefore, from (A.3), (A.8), (A.9), (A.11), (A.12), this proves that problem MQC is an exact quadratically constrained reformulation of the nonlinear problem NLP2.

## 9. ACKNOWLEDGMENTS

The authors gratefully acknowledge the financial support from Johnson and Johnson, the Fulbright Program and the Ministerio de Educación de Argentina, and the Center for Advanced Process Decision-making (CAPD) from Carnegie Mellon University. We would also like to thank Kyle Harshbarger for his useful comments during EWO meetings at CMU and Alev Kaya for her valuable insights on MEIO.